\documentstyle[12pt, amsfonts]{amsart}
\begin{document}
\def\R{{\mathbb R}}
\def\Z{{\mathbb Z}}
\def\C{{\mathbb C}}
\newcommand{\trace}{\rm trace}
\newcommand{\Ex}{{\mathbb{E}}}
\newcommand{\Prob}{{\mathbb{P}}}
\newcommand{\E}{{\cal E}}
\newcommand{\F}{{\cal F}}
\newtheorem{df}{Definition}
\newtheorem{theorem}{Theorem}
\newtheorem{lemma}{Lemma}
\newtheorem{pr}{Proposition}
\newtheorem{co}{Corollary}
\def\n{\nu}
\def\sign{\mbox{ sign }}
\def\a{\alpha}
\def\N{{\mathbb N}}
\def\A{{\cal A}}
\def\L{{\cal L}}
\def\X{{\cal X}}
\def\F{{\cal F}}
\def\c{\bar{c}}
\def\v{\nu}
\def\d{\delta}
\def\diam{\mbox{\rm dim}}
\def\vol{\mbox{\rm Vol}}
\def\b{\beta}
\def\t{\theta}
\def\l{\lambda}
\def\e{\varepsilon}
\def\colon{{:}\;}
\def\pf{\noindent {\bf Proof :  \  }}
\def\endpf{ \begin{flushright}
$ \Box $ \\
\end{flushright}}
%%%%%%%%%%%%%%%%%%%%%%%%%%%%%%%%%%%%%%%%%%%%%%%%%%%%%%%%%%%%%%%%%%%%

\title[Stability of volume comparison]
{Stability of volume comparison for complex convex bodies}

\author{Alexander Koldobsky}

\address{Department of Mathematics\\ 
University of Missouri\\
Columbia, MO 65211}

\email{koldobskiya@@missouri.edu}



%%%%%%%%%%%%%%%%%%%%%%%%%%%%%%%%%%%%%%%%%%%%%%%%%%%%%%%%%%%%%%%%%%%%
\begin{abstract} We prove the stability of the affirmative part of the solution to the complex Busemann-Petty
problem. Namely, if $K$ and $L$ are origin-symmetric convex bodies in $\C^n,$\ $n=2$ or $n=3,$ $\e>0$ and 
 $\vol_{2n-2}(K\cap H) \le \vol_{2n-2}(L\cap H) + \e$ for any complex hyperplane $H$ in $\C^n,$ then 
$\left(\vol_{2n}(K)\right)^{\frac{n-1}n} \le \left(\vol_{2n}(L)\right)^{\frac{n-1}n} + \e,$ where $\vol_{2n}$ is the 
volume in $\C^n,$ which is identified with $\R^{2n}$ in the natural way. 
\end{abstract}  
\maketitle
%%%%%%%%%%%%%%%%%%%%%%%%%%%%%%%%%%%%%%%%%%%%%%%%%%%%%%%%%%%%%%%%%%%%

\section{Introduction}

The Busemann-Petty problem, posed in 1956 (see \cite{BP}),
asks the following question. Suppose that $K$ and $L$ are origin symmetric
convex bodies in $\R^n$ such that
$$\vol_{n-1}(K\cap H) \le \vol_{n-1}(L\cap H)$$
for every hyperplane $H$ in $\R^n$ containing the origin. Does it
follow that
$$\vol_n(K) \le \vol_n(L)?$$
The answer is affirmative if $n\le 4$ and negative if $n\ge 5.$
The solution was completed in the end of the 90's as the result of
a sequence of papers \cite{LR}, \cite{Ba}, \cite{Gi}, \cite{Bo}, 
\cite{L}, \cite{Pa}, \cite{G1}, \cite{G2}, \cite{Z1}, \cite{Z2}, \cite{K1}, \cite{K2}, \cite{Z3},
\cite{GKS} ; see \cite[p. 3]{K3} or \cite[p. 343]{G3} for the history of the solution.

The complex version of the Busemann-Petty problem was solved in \cite{KKZ}, the 
answer is affirmative for convex bodies in $\C^n$ when $n\le 3$, and it is negative
for $n\ge 4.$ To formulate the complex version, we need several definitions.

For $\xi\in \C^n,\ |\xi|=1,$ denote by
$$H_\xi = \{ z\in \C^n:\ (z,\xi)=\sum_{k=1}^n z_k\overline{\xi_k} =0\}$$
the complex hyperplane through the origin perpendicular to $\xi.$

Origin symmetric convex bodies in $\C^n$ are the unit balls of norms on $\C^n.$
We denote by $\|\cdot\|_K$
the norm corresponding to the body $K:$
$$K=\{z\in \C^n:\ \|z\|_K\le 1\}.$$
In order to define volume, we identify $\C^n$ with $\R^{2n}$ using the mapping
$$\xi = (\xi_1,...,\xi_n)=(\xi_{11}+i\xi_{12},...,\xi_{n1}+i\xi_{n2})
 \mapsto  (\xi_{11},\xi_{12},...,\xi_{n1},\xi_{n2}).$$
 Under this mapping the hyperplane $H_\xi$ turns into a $(2n-2)$-dimensional subspace of
 $\R^{2n}.$ 
 
 Since norms on $\C^n$ satisfy the equality
$$\|\lambda z\| = |\lambda|\|z\|,\quad \forall z\in \C^n,\  \forall\lambda \in \C,$$
origin symmetric complex convex bodies correspond to those origin symmetric convex bodies
$K$  in $\R^{2n}$ that are invariant
 with respect to any coordinate-wise two-dimensional rotation, namely for each $\theta\in [0,2\pi]$
 and each $\xi= (\xi_{11},\xi_{12},...,\xi_{n1},\xi_{n2})\in \R^{2n}$
  \begin{equation} \label{rotation}
  \|\xi\|_K =
 \|R_\theta(\xi_{11},\xi_{12}),...,R_\theta(\xi_{n1},\xi_{n2})\|_K,
 \end{equation}
 where $R_\theta$ stands for  the counterclockwise rotation of $\R^2$ by the angle
 $\theta$ with respect to the origin. We shall simply say that $K$ {\it is invariant with respect
 to all $R_\theta$} if it satisfies  (\ref{rotation}).

 The complex Busemann-Petty problem can be formulated as follows:
 suppose $K$ and $L$ are origin symmetric invariant with respect to all $R_\theta$
 convex bodies in $\R^{2n}$ such that
 $$\vol_{2n-2}(K\cap H_\xi)\le \vol_{2n-2}(L\cap H_\xi)$$
 for each $\xi$ from the unit sphere $S^{2n-1}$ of $\R^{2n}.$ Does it follow that
 $$\vol_{2n}(K) \le \vol_{2n}(L) ?$$
 
 As mentioned above, the answer is affirmative if and only if $n\le 3.$ In this article we prove the stability 
 of the affirmative part of the solution:
 
 \begin{theorem} \label{main} Suppose that $\e>0$,  $K$ and $L$ are origin-symmetric
invariant with respect to all $R_\theta$ convex bodies  bodies in $\R^{2n},$ $n=2$ or $n=3.$ 
If for every $\xi\in S^{2n-1}$
\begin{eqnarray}\label{sect1}
\vol_{2n-2}(K\cap H_\xi)\le \vol_{2n-2}(L\cap H_\xi)+\e,
\end{eqnarray}
then
$$\vol_{2n}(K)^{\frac{n-1}n}  \le \vol_{2n}(L)^{\frac{n-1}n} + \e.$$
\end{theorem}
The result does not hold for $n>3,$ simply because the answer to the complex Busemann-Petty problem
in these dimensions is negative; see \cite{KKZ}.

It immediately follows from Theorem \ref{main} that

\begin{co} \label{n3} If $n=2$ or $n=3,$ then for any origin-symmetric invariant with respect to all $R_\theta$
convex bodies $K,L$ in $\R^{2n},$ 
$$\left|\vol_{2n}(K)^{\frac{n-1}n} - \vol_{2n}(L)^{\frac{n-1}n}\right|$$
$$ \le  \max_{\xi\in S^{2n-1}} \left|\vol_{2n-2}(K\cap H_\xi)- \vol_{2n-2}(L\cap H_\xi)\right|.$$
\end{co}

Note that stability in comparison problems for volumes of convex bodies was studied
in \cite{K5}, where it was proved for the original (real) Busemann-Petty problem.

For other results related to the complex Busemann-Petty problem see \cite{R}, \cite{Zy1}, \cite{Zy2}.

\section{Proofs}
We use the techniques of the Fourier approach
to sections of convex bodies; see \cite{K3} and \cite{KY} for details. 

The Fourier transform of a
distribution $f$ is defined by $\langle\hat{f}, \phi\rangle= \langle f, \hat{\phi} \rangle$ for
every test function $\phi$ from the Schwartz space $ \mathcal{S}$ of rapidly decreasing infinitely
differentiable functions on $\R^n$. 

If $K$ is a convex body and $0<p<n,$
then $\|\cdot\|_K^{-p}$  is a locally integrable function on $\R^n$ and represents a distribution. 
Suppose that $K$ is infinitely smooth, i.e. $\|\cdot\|_K\in C^\infty(S^{n-1})$ is an infinitely differentiable 
function on the sphere. Then by \cite[Lemma 3.16]{K3}, the Fourier transform of $\|\cdot\|_K^{-p}$  
is an extension of some function $g\in C^\infty(S^{n-1})$ to a homogeneous function of degree
$-n+p$ on $\R^n.$ When we write $\left(\|\cdot\|_K^{-p}\right)^\wedge(\xi),$ we mean $g(\xi),\ \xi \in S^{n-1}.$
If $K,L$ are infinitely smooth star bodies, the following spherical version of Parseval's
formula was proved in \cite{K4} (see \cite[Lemma 3.22]{K3}):  for any $p\in (-n,0)$
\begin{equation}\label{parseval}
\int_{S^{n-1}} \left(\|\cdot\|_K^{-p}\right)^\wedge(\xi) \left(\|\cdot\|_L^{-n+p}\right)^\wedge(\xi) =
(2\pi)^n \int_{S^{n-1}} \|x\|_K^{-p} \|x\|_L^{-n+p}\ dx.
\end{equation}

A distribution is called {\it positive definite} if its Fourier transform is a positive distribution in
the sense that $\langle \hat{f},\phi \rangle \ge 0$ for every non-negative test function $\phi.$

The Fourier transform formula for the volume of complex hyperplane sections was proved in \cite{KKZ}: 
\begin{pr} \label{volume-ft}  Let $K$ be an infinitely smooth origin symmetric
invariant with respect to $R_\theta$ convex
body in $\R^{2n}, n\ge 2.$  For every $\xi\in S^{2n-1},$ we have
\begin{equation}\label{vft}
\vol_{2n-2}(K\cap H_\xi) = \frac{1}{4\pi(n-1)} \left(\|\cdot\|_K^{-2n+2}\right)^\wedge(\xi).
\end{equation}
\end{pr}

We also use the result of Theorem 3 from \cite{KKZ}. It is formulated in \cite{KKZ} in terms of embedding in $L_{-p},$
which is equivalent to our formulation below. However, the reader does not need to worry about
embeddings in $L_{-p},$ because the proof of Theorem 3 in \cite{KKZ} directly establishes the following:
\begin{pr}\label{Lem:pos-def} Let  $n\ge 3.$ For every origin symmetric invariant with respect to $R_\theta$
convex body $K$ in $\R^{2n},$ the function $\|\cdot\|_K^{-2n+4}$ represents a positive definite distribution.
\end{pr}

Let us formulate precisely what we are going to use later.  The case $n=2$ follows from Proposition \ref{volume-ft} 
(obviously, the volume is positive), the case $n=3$ is immediate from Proposition \ref{Lem:pos-def}.
\begin{co} \label{pos} If $n=2$ or $n=3,$ then for every origin symmetric infinitely smooth invariant with respect to $R_\theta$
convex body $K$ in $\R^{2n},$ $\left(\|\cdot\|_K^{-2}\right)^\wedge$ is a non-negative infinitely smooth function 
on the sphere $S^{2n-1}.$
\end{co}
We need the following simple fact:
\begin{lemma} \label{gammafunction} For every $n\in \N,$
$$\left(\Gamma(n)\right)^{\frac 1n} \le n^{\frac{n-1}n}.$$
\end{lemma}
\pf By log-convexity of the $\Gamma$-function (see \cite[p.30]{K3}),
$$\frac{\log(\Gamma(n+1))-\log(\Gamma(1))}{n} \ge \frac{\log(\Gamma(n))-\log(\Gamma(1))}{n-1},$$ 
so
$$\left(\Gamma(n+1)\right)^{\frac{n-1}n} \ge \Gamma(n).$$
Now note that $\Gamma(n+1)=n\Gamma(n).$
\endpf

The polar formula for the volume of a convex body $K$ in $\R^{2n}$ reads as follows (see \cite[p.16]{K3}):
\begin{equation} \label{polar}
\vol_{2n}(K) = \frac 1{2n}\int_{S^{2n-1}} \|x\|_K^{-2n} dx.
\end{equation}
\bigbreak
We are now ready to prove Theorem \ref{main}. 
\medbreak
{\bf Proof of Theorem \ref{main}.} By the approximation argument of \cite[Th. 3.3.1]{S} (see also \cite{GZ}),
we may assume that the bodies $K$ and $L$ are infinitely smooth. Using \cite[Lemma 3.16]{K3}
we get in this case that the Fourier transforms 
$\left(\|\cdot\|_K^{-2n+2}\right)^\wedge, \left(\|\cdot\|_L^{-2n+2}\right)^\wedge,\left(\|\cdot\|_K^{-2}\right)^\wedge$ 
are the extensions of infinitely differentiable functions on the sphere to homogeneous functions on $\R^{2n}.$ 

By (\ref{vft}), the condition (\ref{sect1}) can be written as
$$ \left(\|\cdot\|_K^{-2n+2}\right)^\wedge(\xi) \le \left(\|\cdot\|_L^{-2n+2}\right)^\wedge(\xi) + 4\pi(n-1)\e$$
for every $\xi\in S^{2n-1}.$ Integrating both sides with respect to a non-negative (by Corollary \ref{pos}) density,
we get
$$\int_{S^{2n-1}}  \left(\|\cdot\|_K^{-2n+2}\right)^\wedge(\xi)\left(\|\cdot\|_K^{-2}\right)^\wedge(\xi) d\xi$$
$$\le \int_{S^{2n-1}}  \left(\|\cdot\|_L^{-2n+2}\right)^\wedge(\xi)\left(\|\cdot\|_K^{-2}\right)^\wedge(\xi) d\xi$$
$$+\ 4\pi(n-1)\e  \int_{S^{2n-1}}  \left(\|\cdot\|_K^{-2}\right)^\wedge(\xi) d\xi.$$
By the Parseval formula (\ref{parseval}) applied twice, 
$$(2\pi)^n \int_{S^{2n-1}}  \|x\|_K^{-2n} dx \le (2\pi)^n \int_{S^{2n-1}}  \|x\|_L^{-2n+2}\|x\|_K^{-2} dx$$  
$$+ \ 4\pi(n-1)\e  \int_{S^{2n-1}}  \left(\|\cdot\|_K^{-2}\right)^\wedge(\xi) d\xi.$$
Estimating the first summand in the right-hand side of the latter inequality by H\"older's inequality,
$$(2\pi)^n \int_{S^{2n-1}}  \|x\|_K^{-2n} dx \le (2\pi)^n \left(\int_{S^{2n-1}} \|x\|_L^{-2n} dx\right)^{\frac {n-1}n}
\left(\int_{S^{2n-1}} \|x\|_K^{-2n} dx\right)^{\frac {1}n} $$
$$+ \ 4\pi(n-1)\e  \int_{S^{2n-1}}  \left(\|\cdot\|_K^{-2}\right)^\wedge(\xi) d\xi.$$
and using the polar formula for the volume (\ref{polar}),
$$(2\pi)^n (2n) \vol_{2n}(K) \le (2\pi)^n (2n) \left(\vol_{2n}(L)\right)^{\frac {n-1}n}\left(\vol_{2n}(K)\right)^{\frac {1}n}$$
\begin{equation}\label{eq11}
+  \ 4\pi(n-1)\e  \int_{S^{2n-1}}  \left(\|\cdot\|_K^{-2}\right)^\wedge(\xi) d\xi.
\end{equation}
We now estimate the second summand in the right-hand side. First we use the formula
for the Fourier transform  (in the sense of distributions; see \cite[p.194]{GS})
$$\left(| \cdot |_2^{-2n+2}\right)^\wedge(\xi) = \frac{4\pi^{n}}{\Gamma(n-1)},$$
where $|\cdot |_2$ is the Euclidean norm in $\R^{2n}$ and $\xi \in S^{2n-1}.$ We get
$$4\pi(n-1)\e  \int_{S^{2n-1}}  \left(\|\cdot\|_K^{-2}\right)^\wedge(\xi) d\xi $$
$$= \frac{4\pi(n-1)\Gamma(n-1)\e}{4\pi^n}  \int_{S^{2n-1}}  \left(\|\cdot\|_K^{-2}\right)^\wedge(\xi) 
\left(| \cdot |_2^{-2n+2}\right)^\wedge(\xi)d\xi,$$
and by Parseval's formula (\ref{parseval}) and H\"older's inequality,
$$= \frac{(2\pi)^n\e \Gamma(n)}{\pi^{n-1}} \int_{S^{2n-1}}  \|x\|_K^{-2} dx$$
$$\le \frac{(2\pi)^n\e \Gamma(n)}{\pi^{n-1}}\left(\int_{S^{2n-1}} \|x\|_K^{-2n} dx \right)^{\frac 1n} \left| S^{2n-1}\right|^{\frac{n-1}n},$$
where $\left|S^{2n-1}\right| = (2\pi^n)/\Gamma(n)$ is the surface area of the unit sphere in $\R^{2n}.$
By the polar formula for the volume, the latter is equal to
$$(2\pi)^n (2n) \e \left(\vol_{2n}(K)\right)^{\frac{1}n}\frac{\left(\Gamma(n)\right)^{\frac 1n}}{n^{\frac{n-1}n}}
\le (2\pi)^n (2n) \e \left(\vol_{2n}(K)\right)^{\frac{1}n}$$
by Lemma \ref{gammafunction}.
Combining this with (\ref{eq11}), we get the result. \qed
\bigbreak
We finish with the following ``separation" property (see \cite{K5} for more results of this kind). 
Note that for any $x\in S^{2n-1},$ $\|x\|_K^{-1}=\rho_K(x)$
is the radius of $K$ in the direction $x,$ and denote by
$$r(K) = \frac{\min_{x\in S^{2n-1}} \rho_K(x)}{\left(\vol_{2n}(K)\right)^{\frac 1{2n}}}$$
the normalized inradius of $K.$ Clearly,  for every $x\in S^{2n-1}$ we have 
$$\|x\|_K^{-1}\ge r(K)\left(\vol_{2n}(K)\right)^{\frac 1{2n}}.$$
 \begin{theorem} \label{separ} Suppose that $\e>0$,  $K$ and $L$ are origin-symmetric
invariant with respect to all $R_\theta$ convex bodies  bodies in $\R^{2n},$ $n=2$ or $n=3.$ 
If for every $\xi\in S^{2n-1}$
$$
\vol_{2n-2}(K\cap H_\xi)\le \vol_{2n-2}(L\cap H_\xi)-\e,
$$
then
$$\vol_{2n}(K)^{\frac{n-1}n}  \le \vol_{2n}(L)^{\frac{n-1}n} - \frac{\pi r^2(K)}{n}\e.$$
\end{theorem}

\pf We follow the lines of the proof of Theorem \ref{main} to get 
$$(2\pi)^n (2n) \vol_{2n}(K) \le (2\pi)^n (2n) \left(\vol_{2n}(L)\right)^{\frac {n-1}n}\left(\vol_{2n}(K)\right)^{\frac {1}n}$$
\begin{equation}\label{eq21}
-  \ 4\pi(n-1)\e  \int_{S^{2n-1}}  \left(\|\cdot\|_K^{-2}\right)^\wedge(\xi) d\xi.
\end{equation}
We now need a lower estimate for
$$4\pi(n-1)\e  \int_{S^{2n-1}}  \left(\|\cdot\|_K^{-2}\right)^\wedge(\xi) d\xi.$$
Similarly to how it was done in Theorem \ref{main}, we write the latter as
$$\frac{(2\pi)^n\e \Gamma(n)}{\pi^{n-1}} \int_{S^{2n-1}}  \|x\|_K^{-2} dx \ge
\frac{(2\pi)^n\e \Gamma(n)r^2(K)\left(\vol_{2n}(K)\right)^{\frac 1n}}{\pi^{n-1}}\left|S^{2n-1}\right|. \qed$$

\bigbreak

{\bf Acknowledgement.} The author wishes to thank
the US National Science Foundation for support through 
grants DMS-0652571 and DMS-1001234.


\begin{thebibliography}{99}

\bibitem[Ba]{Ba} { K.~Ball}, \textit{Some remarks on the geometry of convex sets},
Geometric aspects of functional analysis (1986/87), Lecture Notes in Math. \textbf{1317}, 
Springer-Verlag, Berlin-Heidelberg-New York, 1988, 224--231.

\bibitem[Bo]{Bo} { J.~Bourgain}, \textit{On the Busemann-Petty problem for
perturbations of the ball}, Geom. Funct. Anal. \textbf{ 1} (1991), 1--13.

\bibitem[BP]{BP} { H.~Busemann and C.~M.~Petty}, \textit{Problems on convex bodies},
Math. Scand. \textbf{4} (1956), 88--94.

\bibitem[G1]{G1} { R.~J.~Gardner}, \textit{Intersection bodies and the 
Busemann-Petty problem}, Trans. Amer. Math. Soc. \textbf{342} (1994), 
435--445.

\bibitem[G2]{G2} { R.~J.~Gardner}, \textit{A positive answer to the Busemann-Petty 
problem in three dimensions}, Annals of Math. \textbf{140} (1994), 435--447.

\bibitem[G3]{G3} { R.~J.~Gardner}, \textit{Geometric tomography}, Second edition,
Cambridge University Press, Cambridge, 2006. 

\bibitem[GKS]{GKS} { R.~J.~Gardner, A.~Koldobsky and Th.~Schlumprecht},
\textit{An analytic solution to the Busemann-Petty 
problem on sections of convex bodies}, Annals of Math. \textbf{149} (1999),
691--703.

\bibitem[GS]{GS} { I.~M.~Gelfand and G.~E.~Shilov}, \textit{Generalized functions, vol.~1. 
Properties and operations}, Academic Press, New York, 1964.

\bibitem[Gi]{Gi} { A.~Giannopoulos}, \textit{A note on a problem of H.~Busemann and 
C.~M.~Petty concerning sections of symmetric convex bodies},
Mathematika \textbf{37} (1990),  239--244.

\bibitem[GZ]{GZ} { E.~Grinberg and Gaoyong Zhang}, \textit{Convolutions, 
transforms, and convex bodies}, Proc. London Math. Soc. (3) \textbf{78}  (1999), 77--115.

\bibitem[K1]{K1} { A. Koldobsky}, \textit{Intersection bodies, positive definite 
distributions and the Busemann-Petty problem}, Amer. J. Math. \textbf{120} (1998), 827--840.

\bibitem[K2]{K2} { A. Koldobsky}, \textit{Intersection bodies in ${\R^4}$}, 
Adv. Math. \textbf{136} (1998), 1--14.

\bibitem[K3]{K3} {A.~Koldobsky}, \textit{Fourier analysis in convex geometry},
Amer. Math. Soc., Providence RI, 2005.

\bibitem[K4]{K4} { A.~Koldobsky}, \textit{A generalization of the Busemann-Petty problem 
on sections of convex bodies}, Israel J. Math. 
\textbf{110} (1999), 75--91.

\bibitem[K5]{K5} { A.~Koldobsky}, \textit{Stability in the Busemann-Petty and Shephard  problems}, 
preprint.

\bibitem[KKZ]{KKZ} {A.~Koldobsky, H.~K\"onig and M.~Zymonopoulou}, \textit{The complex Busemann-Petty
problem on sections of convex bodies}, Adv. Math.  \textbf{218} (2008), 352--367.

\bibitem[KY]{KY} {A.~Koldobsky and V.~Yaskin}, \textit{The interface between convex geometry and harmonic analysis}, 
CBMS Regional Conference Series in Mathematics, 108, American Mathematical Society, Providence, RI, 2008.

\bibitem[LR]{LR} { D.~G.~Larman and C.~A.~Rogers}, \textit{The existence of a centrally 
symmetric convex body with central sections that are unexpectedly small},
Mathematika  \textbf{22} (1975), 164--175.  

\bibitem[L]{L} { E.~Lutwak}, \textit{Intersection bodies and dual mixed volumes},
Adv. Math. \textbf{71} (1988), 232--261.

\bibitem[Pa]{Pa} { M.~Papadimitrakis}, \textit{On the Busemann-Petty problem
about convex, centrally symmetric bodies in $\R^n$},
Mathematika \textbf{39} (1992), 258--266.

\bibitem[R]{R} {B.~Rubin},  \textit{Comparison of volumes of convex bodies in real, complex, and quaternionic spaces},
 Adv. Math. \textbf{225} (2010), 1461--1498.

\bibitem[S]{S} { R.~Schneider}, \textit{Convex bodies: the Brunn-Minkowski theory},
Cambridge University Press, Cambridge, 1993.
 
 \bibitem[Z1]{Z1} { Gaoyong Zhang}, \textit{Centered bodies and dual mixed volumes},
Trans. Amer. Math. Soc. \textbf{345} (1994), 777--801. 

\bibitem[Z2]{Z2} { Gaoyong Zhang}, \textit{Intersection bodies and Busemann-Petty
inequalities in $\R^4$}, Annals of Math. \textbf{140} (1994), 331--346.

\bibitem[Z3]{Z3} { Gaoyong Zhang}, \textit{A positive answer to the Busemann-Petty 
problem in four dimensions}, Annals of Math. \textbf{149} (1999), 535--543.

\bibitem[Zy1]{Zy1} {M.~Zymonopoulou}, \textit{The modified complex Busemann-Petty problem on sections of convex bodies}, 
Positivity \textbf{13} (2009), no. 4, 717--733.
 
\bibitem[Zy2]{Zy2}{M.~Zymonopoulou}, \textit{The complex Busemann-Petty problem for arbitrary measures},
 Arch. Math. (Basel) \textbf{91} (2008), no. 5, 436--449.


\end{thebibliography}
\end{document}